\documentclass[11pt,renqo,reqno]{amsart}

\usepackage{graphicx} 
\usepackage{verbatim} 
\usepackage{amsmath, amsthm, amssymb, amsrefs} 
\usepackage{anysize} 
\usepackage[autostyle,english=american]{csquotes} 
\usepackage[shortlabels]{enumitem} 
\usepackage{hyperref} 
\usepackage{bm} 
\usepackage{color} 
\usepackage[parfill]{parskip} 

\marginsize{2.25cm}{2.25cm}{2.25cm}{2.25cm} 
\MakeOuterQuote{"} 
\setlist[enumerate]{topsep=0pt,itemsep=0ex,partopsep=1ex,parsep=1ex,leftmargin=*} 
\numberwithin{equation}{section} 

\newtheorem{theorem}{Theorem}

\theoremstyle{remark}

\theoremstyle{definition}

\def\to{\rightarrow}

\def\Z{\mathbb{Z}}

\def\1{1\!\!1}

\title{Optimally Packing a Large Square by Unit Squares}
\author{Rory McClenagan}
\address{Department of Mathematics and Statistics \\
        University of Northern British Columbia \\
        Prince George, BC V2N4Z9 \\
        Canada}
\date{\today}

\begin{document}
\maketitle
\begin{abstract}
We show that a large square of sidelength $x$ can be packed by unit squares in a manner so that the wasted space $W(x) = O(x^{3/5})$.
\end{abstract}

Let $S(x)$ denote a square of sidelength $x$ for some large $x$. Pack $S(x)$ as efficiently as possible by squares of unit sidelength with disjoint interiors. Let $W(x)$ denote the minimum amount of area left uncovered in $S(x)$ by any such packing. Of course, if $x$ is an integer, then $W(x) =0$. In general, when $x$ is not an integer, orienting the unit squares so that they are parallel to the sides of $S(x)$ in a naive manner will generate a wasted space $W(x)  = O(x \{x\})$. Here, $\{x\}$ is the fractional part of $x$, and if it is bounded away from $0$, the wasted space $W(x)$ simply becomes $O(x)$.

\begin{figure}[h]
    \includegraphics[width=.6\textwidth]{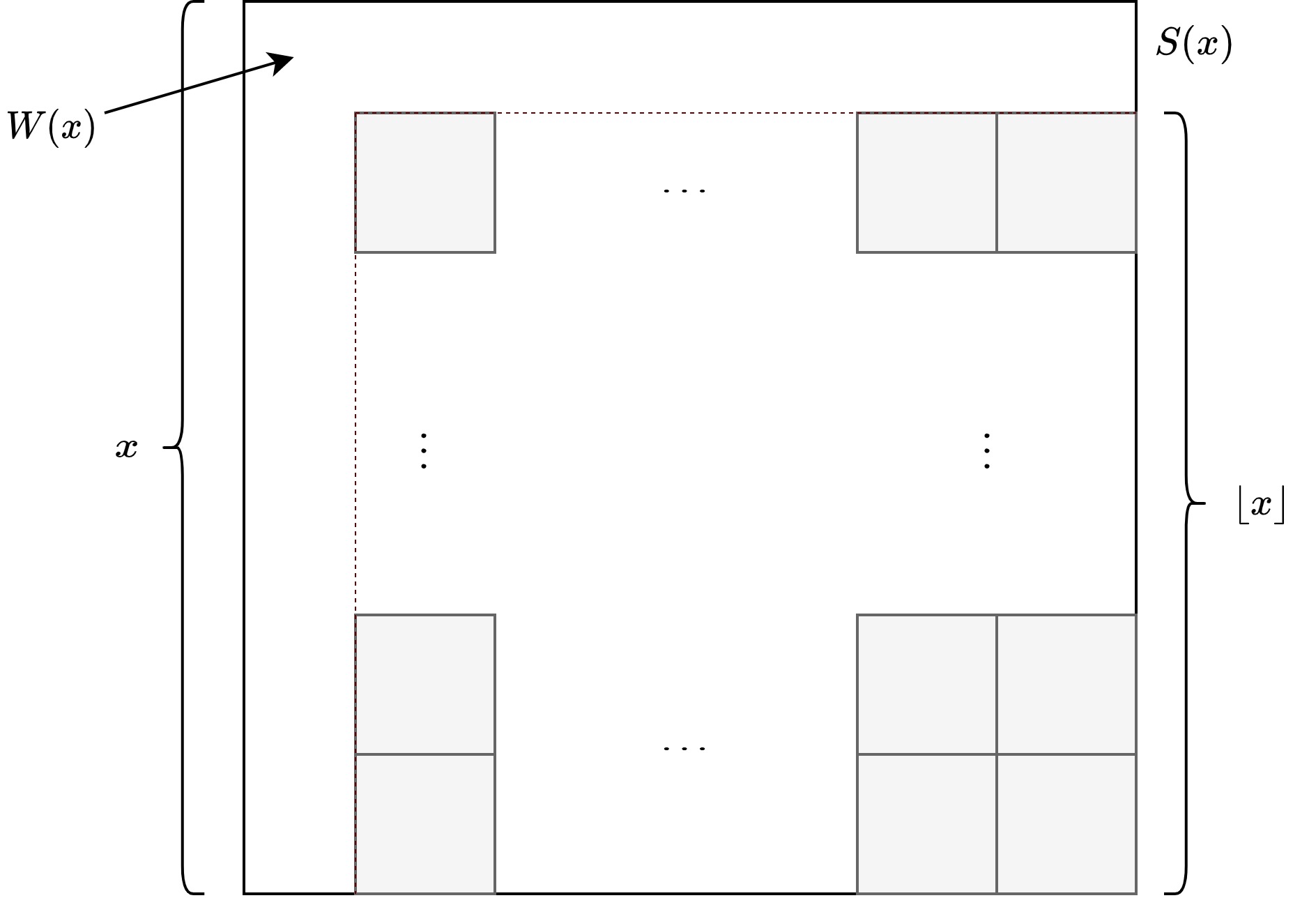}
    \centering
    \caption{The trivial packing of $S(x)$ by squares of unit sidelength. The wasted space $W(x)$ is only $O(x)$ if $\{x\} = x - \lfloor x \rfloor$ is bounded away from $0$.}
    \label{fig-parallel_pack}
\end{figure}

However, if one packs the squares at slight angles, the wasted space can be decreased. In 1975, Paul Erd\"{o}s and Ronald Graham showed  that $W(x) = O(x^{\frac{7}{11}})$ (see \cite{erdos-square_packing}). This prompted discussion around the proper order of growth of $W(x)$. We do know that the bound can at least be no better than $O(x^{1/2})$. This result is due to Roth and Vaughan in 1978 (see \cite{roth-omega_square_packing}) who showed that
    \[W(x) > 10^{-100} \sqrt{x - \lfloor x \rfloor}.\]
On the other end, Montgomery improved the upper bound to $W(x) = O(x^{\frac{3 - \sqrt{3}}{2}})$, according to personal communication (see, for example, \cite{chung_graham-3_5}). In 2009, Fan Chung and Ronald Graham improved it further to $W(x) = O(x^{\frac{3+\sqrt{2}}{7}} \log x)$. Now, recently in 2020, Chung and Graham claimed that this could be improved to $W(x) = O(x^{3/5})$ (see \cite{chung_graham-3_5}). Unfortunately, this result has an error in it, which brings the best known bound back to $W(x) = O(x^{\frac{3+\sqrt{2}}{7}} \log x)$. In this paper, we show that the bound claimed by Chung and Graham in \cite{chung_graham-3_5} is, in fact, correct, using a new algorithm:

\begin{theorem}\label{thm-main}
The wasted space in packing the square $S(x)$ by unit squares is bounded by
    \[W(x) = O(x^{3/5}).\]
\end{theorem}

We will prove Theorem \ref{thm-main} in several stages. First, we pack all of $S(x)$ using "stacks" of unit squares except for a finite number of trapezoidal regions $T$ in Section \ref{s-S_packing}. We then pack such a generic trapezoidal region using two different "sub-algorithms", described in in sections \ref{s-first_packing} and \ref{s-second_packing}. In section \ref{s-proof}, we show how these can be combined to form a packing of $T$. We then optimize our parameters, and demonstrate that the wasted space generated by our packing is $O(x^{3/5})$.

\subsection*{Acknowledgements}
We would like to thank Dr.~Alia Hamieh for helpful comments and proof-reading this paper.

\subsection*{Notation}\label{s-lemmas}
Throughout this paper, we will use the standard asymptotic notation $X= O(Y)$, $X \ll Y$, and $Y \gg X$ to refer to the relation $X \leq C |Y|$. If, instead, we use the notation $X = O_M(Y)$ or $X \ll_M Y$, then the corresponding constant $C$ is allowed to depend on the parameter $M$. We use $X \asymp Y$ if $X \ll Y$ and $Y\ll X$. Finally, we use the notation $X = o(Y)$ if $X/ Y \to 0$ with respect to some explicit or implicit limiting behaviour defined in context.

\section{Packing \texorpdfstring{$S(x)$}{S(x)} Using Stacks of Unit Squares}\label{s-S_packing}
\begin{figure}[t]
    \includegraphics[width=.7\textwidth]{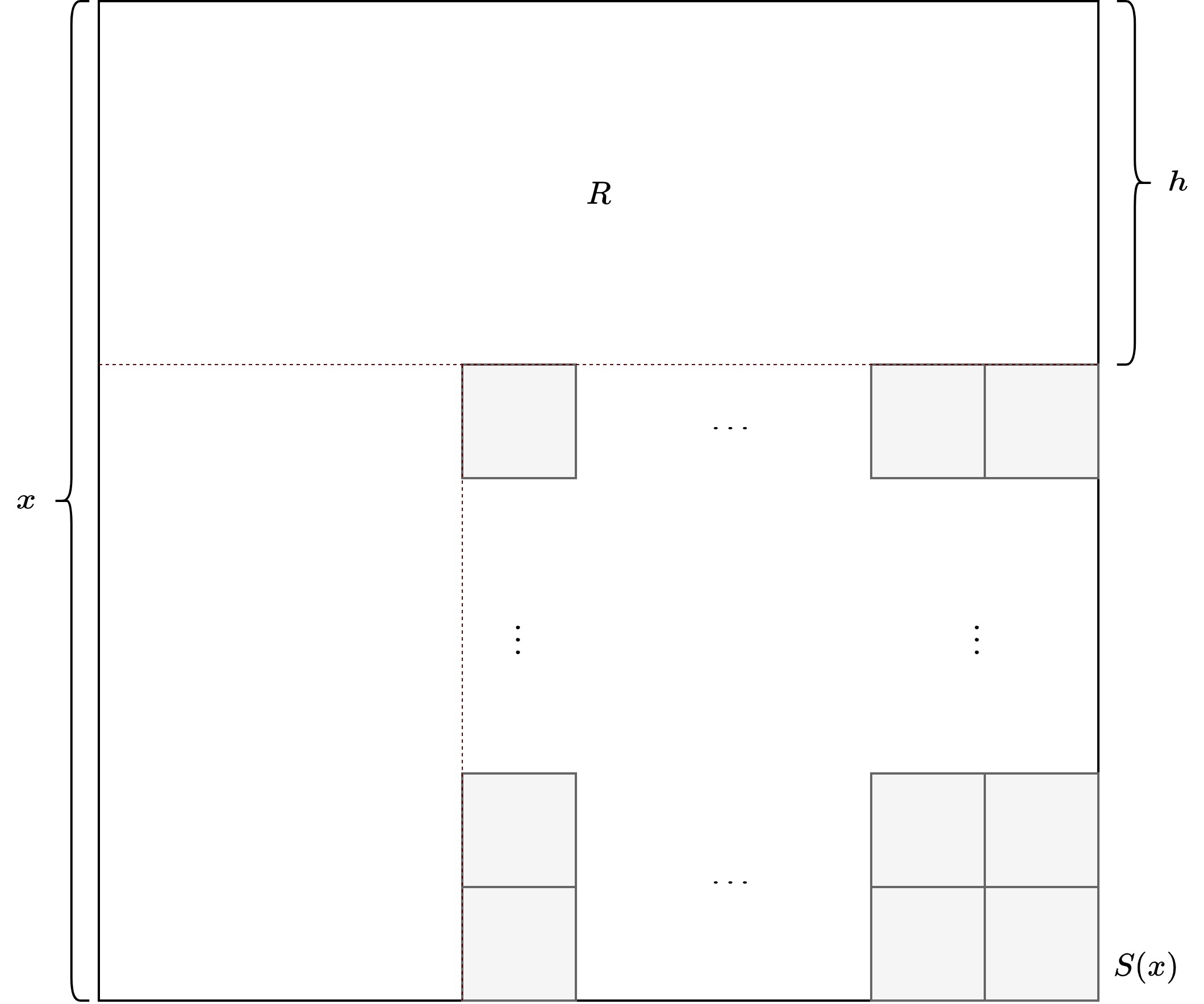}
    \centering
    \caption{We begin by packing $S(x)$ trivially except for two rectangles of width $h$.}
    \label{fig-s_packing}
\end{figure}

We begin by packing $S(x)$ in a trivial manner by placing unit squares snugly with sides parallel to $S(x)$ starting at one corner until the only unpacked region is two rectangles with width $h$ (see Figure \ref{fig-s_packing}). Note that we can fix $h$ to the nearest $O(1)$ by changing the number of squares that we pack during this stage, which is something we will do at the end of the proof.

\begin{figure}[t]
    \includegraphics[width=\textwidth]{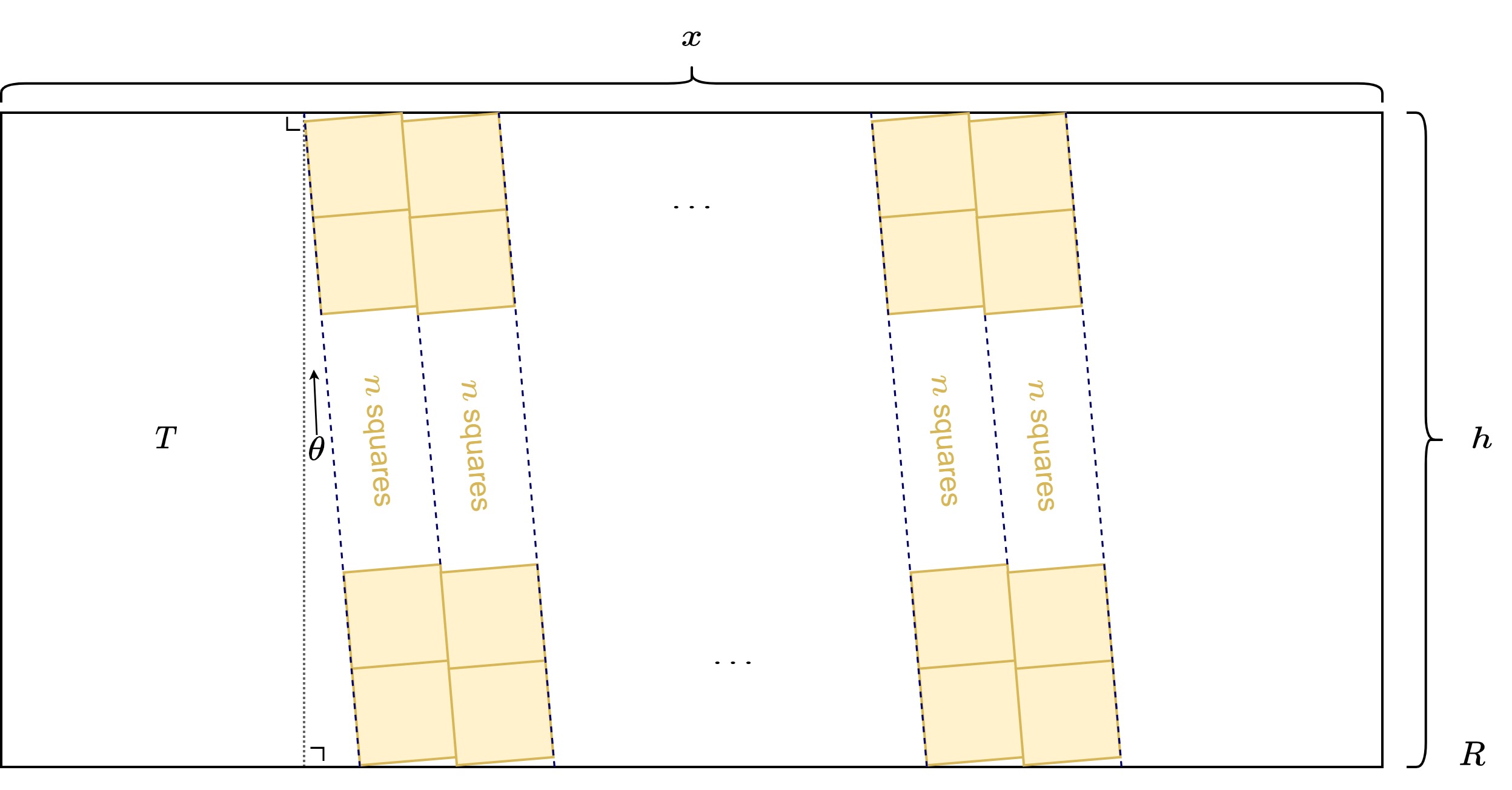}
    \centering
    \caption{We pack $R$ by $n$-length parallel stacks of unit squares inclined at an angle $\theta$ such that each stack touches both the top and the bottom of $R$. We call one of the unpacked trapezoids formed $T$.}
    \label{fig-r_packing}
\end{figure}

Without loss of generality, we will pack only one of the rectangles, which we call $R$. We will assume that the short $h$-length side of $R$ is aligned with the vertical coordinate axis. We then pack $R$ by near-vertical stacks of unit squares of length $n$ inclined so that they touch both sides of $R$ (see Figure \ref{fig-r_packing}). Choose $n \in \Z^+$ so that $1 < n - h \ll 1$. If $\theta$ is the angle of inclination of these stacks, then observe that
    \[\sec(\theta) = \frac{n + \tan \theta}{h} .\]

Taylor series expansion then gives
    \[\theta^2 = \frac{2(n-h)}{h} + O\left ( \frac{\theta}{h} + \theta^4 \right)  .\]
Thus,
\begin{equation}\label{eq-theta_approx}
    \theta \asymp \frac{1}{\sqrt{h}}.
\end{equation}
Pack all of $R$ in such a manner except for a narrow trapezoid on each end of height $h$. Without loss of generality, we will pack only one such trapezoid, which we call $T$. Note that we can fix the width  of $T$ (namely, the width of its smaller side), which we call $w$, to the nearest $O(1)$ by changing the number of $n$-length stacks that we pack during this stage, and will do so at a later time. We will assume that the $h$-length side of $T$ is aligned with the vertical coordinates axis. We will call this the \textit{vertical wall} of $T$. We will call the opposing side the \textit{inclined wall} of $T$. 

Except for a portion that we pack trivially at the top and bottom of $T$, we will pack $T$ in a collection of near horizontal strips that will each have a height of about $2 \theta^{-1}$. Each such strip will be packed using two packing algorithms. Such packings will be replicated all the way down the length of $T$. When we describe these two packing algorithm in sections \ref{s-first_packing} and \ref{s-second_packing}, we will describe them abstractly and not yet fix where we are vertically within $T$, but instead will only refer to the two sides of $T$ as the vertical and inclined walls. We will combine these two packing algorithms together formally in section \ref{s-proof}.

\section{The First Packing Algorithm}\label{s-first_packing}
\subsection*{Setup}
\begin{figure}[t]
    \includegraphics[width=.8\textwidth]{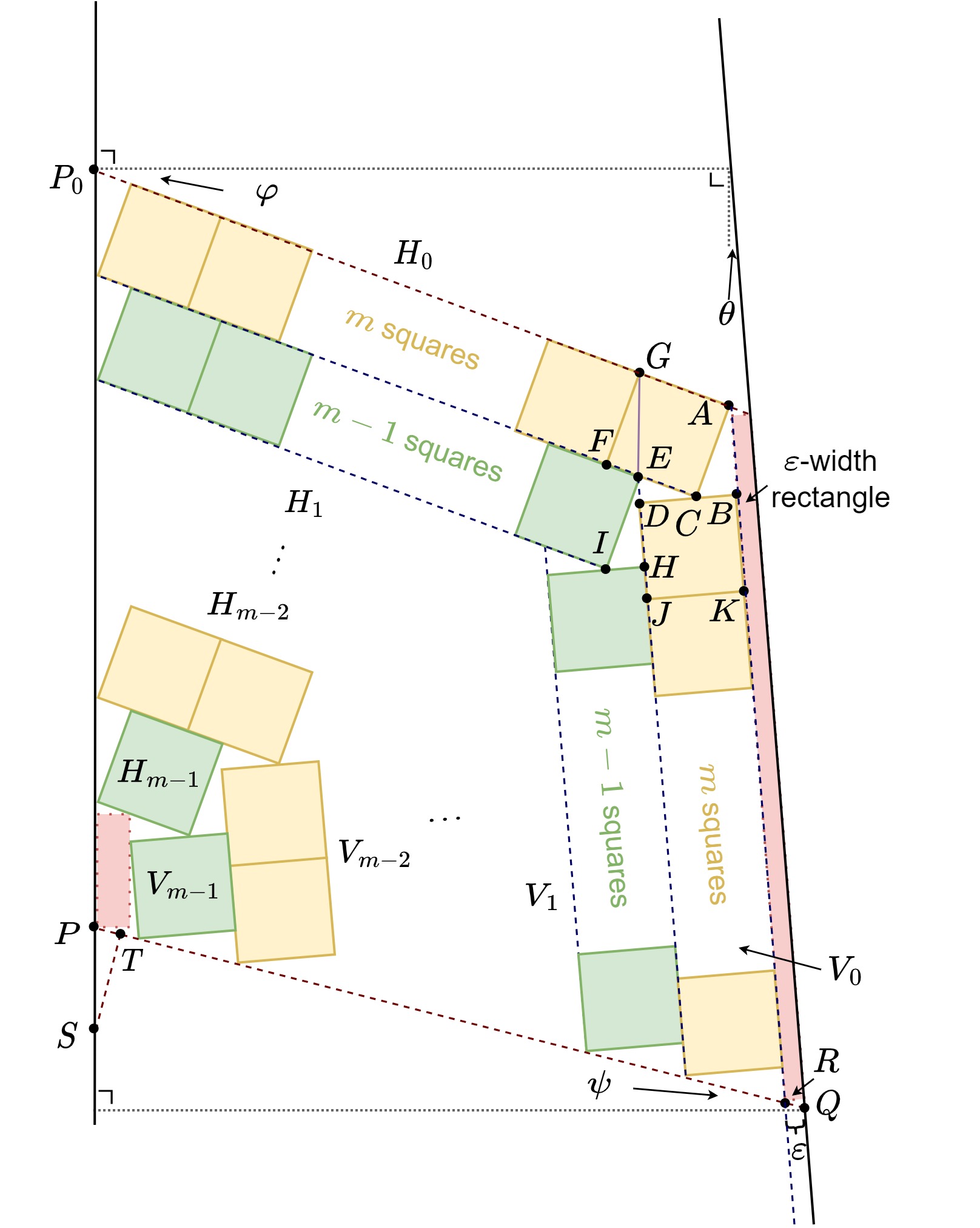}
    \centering
    \caption{The first packing algorithm.}
    \label{fig-packing_alg_1}
\end{figure}

First, place a near-horizontal stack $H_0$ inclined at an angle of $\varphi$ containing $m$ unit squares so that its left edge touches the vertical wall and its right edge is some distance $\varepsilon \geq 0$ from the inclined wall (see Figure \ref{fig-packing_alg_1}). We now place a small rectangle of width $\varepsilon$ along the inclined wall stretching over the entire vertical region that will be covered during the first packing algorithm. This region will not be packed. Place a vertical stack $V_0$ inclined at an angle of $\theta$ with $m$ squares so that it is sung against the $\varepsilon$-width rectangle and parallel to the inclined wall. 

Next place a horizontal stack $H_1$ snugly against $H_0$ containing $m-1$ squares such that it touches the vertical wall. Similarly, place a vertical stack $V_1$ containing $m-1$ squares snugly against $V_0$ and touching $H_1$. We then continue in a similar manner to place the stacks $H_i$ and $V_i$ containing $m-i$ squares for $i=0, \dots, m-1$.

Label the points on the first two horizontal and vertical stacks as in Figure \ref{fig-packing_alg_1}. Construct the line $GE$, noting that this must be vertical by construction. We are going to fix $\varphi$ in such a way that the angle $\angle EDC$ formed is a right-angle. This means that there is a line formed from $E$ to the bottom left corner of $V_0$ inclined at angle $\theta$, just like the inclined wall, which allows us to iteratively pack each consecutive $H_i$ and $V_i$ in the same manner as $H_0$ and $V_0$.

\subsection*{Estimating \texorpdfstring{$\varphi$}{phi}}
Observe that
    \[\angle FGE = \varphi, \quad \angle ECD = \varphi + \theta, \quad \text{and} \quad \angle CAB = \varphi + \theta,\]
and so
\begin{equation}\label{eq-EC}
    EC = 1 - \tan \varphi \quad \text{and} \quad DC = 1 - \sin(\varphi + \theta).
\end{equation}
Thus, $\angle EDC$ is a right angle if and only if
\begin{equation}\label{eq-varphi}
    \left(1 - \tan \varphi \right) \cos (\varphi + \theta) = 1 - \sin (\varphi + \theta).
\end{equation}
Fix $\varphi$ to be the solution of (\ref{eq-varphi}). We now want to determine an approximation for $\varphi$ in terms of $\theta$. Applying Taylor series expansions gives
    \[ \left( 1- \varphi + O(\varphi^3) \right) \left( 1 - \frac{\varphi^2}{2} - \varphi^2\theta^2 - \frac{\theta^2}{2} + O\left((\varphi + \theta)^4\right) \right) = 1 - \varphi - \theta + O \left((\varphi + \theta)^3\right).\]
Aggregating lower order terms gives
    \[\theta-\frac{\varphi^2}{2} - \varphi^2\theta^2 - \frac{\theta^2}{2} + \frac{\varphi\theta^2}{2} = O \left((\varphi + \theta)^3\right).\]
This implies that $\varphi^2 + O(\varphi^3) = 2\theta + O(\theta^2)$. It follows that $\varphi^2 \asymp \theta$ allowing us to write $\varphi = \sqrt{2\theta + O(\theta^{3/2})}$. Taylor series expansion gives
\begin{equation}\label{eq-varphi_approx}
    \varphi = \sqrt{2\theta} + O(\theta) \quad \text{and} \quad \theta = \frac{1}{2} \varphi^2 + O(\varphi^3).
\end{equation}

\subsection*{Estimating \texorpdfstring{$\psi$}{psi}}
Now, we have fixed $\varphi$ (and estimated it in terms of $\theta$) such that the packing setup shown in Figure \ref{fig-packing_alg_1} is feasible and can be iterated for each $H_i$ and $V_i$. This automatically fixes the angle $\psi$ formed at the bottom of this region that we are packing (see Figure \ref{fig-packing_alg_1}). We now want to estimate this angle. Observe that $\angle HKJ = \psi + \theta$, so
\begin{equation}\label{eq-DH}
    DH = 1 - \tan (\psi + \theta).
\end{equation}
On the other hand since $IEH$ and $ECD$ are similar, so
    \[ED + DH = \frac{DC}{EC}. \]
Since $\angle ECD = \varphi + \theta$, we can rearrange this as
    \[DH = \cos(\varphi + \theta) - EC \sin (\varphi + \theta).\]
Combining this with (\ref{eq-EC}) gives
    \[DH = \cos ( \varphi + \theta)  - (1 - \tan \varphi) \sin (\varphi + \theta)\]
Comparing this with (\ref{eq-DH}) allows us to see that $\psi$ is the unique solution of
\begin{equation}\label{eq-psi}
    1 - \tan (\psi + \theta) = \cos ( \varphi + \theta)  - (1 - \tan \varphi) \sin (\varphi + \theta).
\end{equation}
We can Taylor expand both sides of this equation, simplifying our error term using (\ref{eq-varphi_approx}):
    \[1 - \psi - \theta= 1 - \frac{1}{2} \varphi^2 - (\varphi + \theta) + \varphi^2 + O(\varphi^3 + \psi^3).\]
Solving for $\psi$ gives
    \[\psi = \varphi - \frac{1}{2}\varphi^2  + O(\varphi^3 + \psi^3).\]
It follows that $\psi + O(\psi^3) = \varphi + O(\varphi^2)$ and so $\psi \asymp \varphi$, allowing us to simplify the error term. Combining this observation with (\ref{eq-varphi_approx}), we get
\begin{equation}\label{eq-psi_approx}
    \psi = \varphi - \theta + O(\varphi^3).
\end{equation}
Thus, we have "lost" about $\theta$ in the inclination of the top portion of the region we are packing compared to the bottom.

Now, at this stage, the only quantity we have fixed is $\varphi$, which was chosen in such a way so that we could perform our first packing algorithm. In the second packing algorithm we will want to start by packing $m+1$ square along the length $PR$ compared to the $m$ squares we packed at the beginning of the first packing algorithm (note that $R$ is chosen to be the point on $\overline{PQ}$ such that the distance from $R$ to the inclined wall is $\varepsilon$). To do this, we need to choose $m$ large enough so that the first packing algorithm takes up sufficient vertical space, thus and forces $PR$ to be able to fit $m+1$ squares along it. 

\subsection*{Estimating \texorpdfstring{$m$}{m}}
We want to choose $m$ such that $\overline{PR}$ is just long enough to fit a near-horizontal stack of $m+1$ squares inclined at an angle of $\psi$. Choose points $S$ on the vertical wall and $T$ on $\overline{PR}$ such that $ST = 1$ and $\angle PTS$ is a right-angle (equivalently, $\angle TSP = \psi$). Define $m$ to be the least integer such that $TR \geq m+1$. Note that the distance between the inclined and the vertical walls increases at a rate of $O(\theta)$. Thus, increasing $m$ by $1$ will increase $TR$ by $O(\theta)$. Thus, our definition of $m$ automatically gives 
    \[TR = m + 1 + O(\theta).\] 
It is important to note that this definition of $m$ (and the length $TR$) is \textit{independent} of $\varepsilon$, although of course $TQ$ depends on $\varepsilon$:
\begin{equation}\label{eq-tq_approx}
    TQ = m+1 + O(\theta + \varepsilon).
\end{equation}
We now want to find an approximate expression for $m$ in terms of our known quantities. Observe that the horizontal distance from the vertical wall to the lower-left corner of $V_{m-1}$ is the same as the horizontal distance from $F$ to $J$, which is $O(\varphi)$. Thus, the distance from $P$ to the lower-left corner of $V_{m-1}$ is also $O( \varphi)$, we have:
    \[TR = \frac{m}{\cos(\psi +\theta)} + O(\varphi)\]
However, once again, we know that the difference between $TR$ and $m+1$ is at most $O(\theta)$. Thus,
    \[\frac{m}{\cos (\psi + \theta)} - m = 1 + O(\varphi).\]
Upon Taylor expansion, the left-hand side becomes $\frac{1}{2} \psi^2 m + O(\psi \theta m)$. Applying (\ref{eq-psi_approx}) and (\ref{eq-varphi_approx}) thus gives
    \[\theta m + O(\varphi^3 m) = 1 + O(\varphi).\]
Clearly this implies that $m \asymp \theta^{-1}$, which allows us to simplify the error term, giving us an approximate expression for $m$:
\begin{equation}\label{eq-m_approx}
    m = \theta^{-1} + O(\varphi^{-1}).
\end{equation}

\subsection*{Summary}
This completes the analysis of the first packing stage. We make two important observations. First, observe that the quantities $\varphi, \psi$ and $m$ are all fixed based solely on $\theta$, independent of $\varepsilon$. Second, note that the angle $\varphi$ at the top of the packed squares differs from the angle $\psi$ at the bottom of the packed squares by about $\theta$ (see \ref{eq-psi_approx}). This is the reason we cannot simply iterate this first packing algorithm all the way down the trapezoid $T$. Indeed, if we were to reset and apply the first packing algorithm successively each application would generate a triangle of height $O(m)$ and angle $\varphi - \psi \sim \theta$ having area $\theta^{-1}$. This is an order of magnitude larger than the area of the other wasted space that is generated during the algorithm we just described. We instead have to pair this first algorithm with another algorithm that "undoes" the angle change of $\theta$. We do this in the next section.

\section{The Second Packing Algorithm}\label{s-second_packing}
\subsection*{Setup}
\begin{figure}[t]
    \includegraphics[width=.8\textwidth]{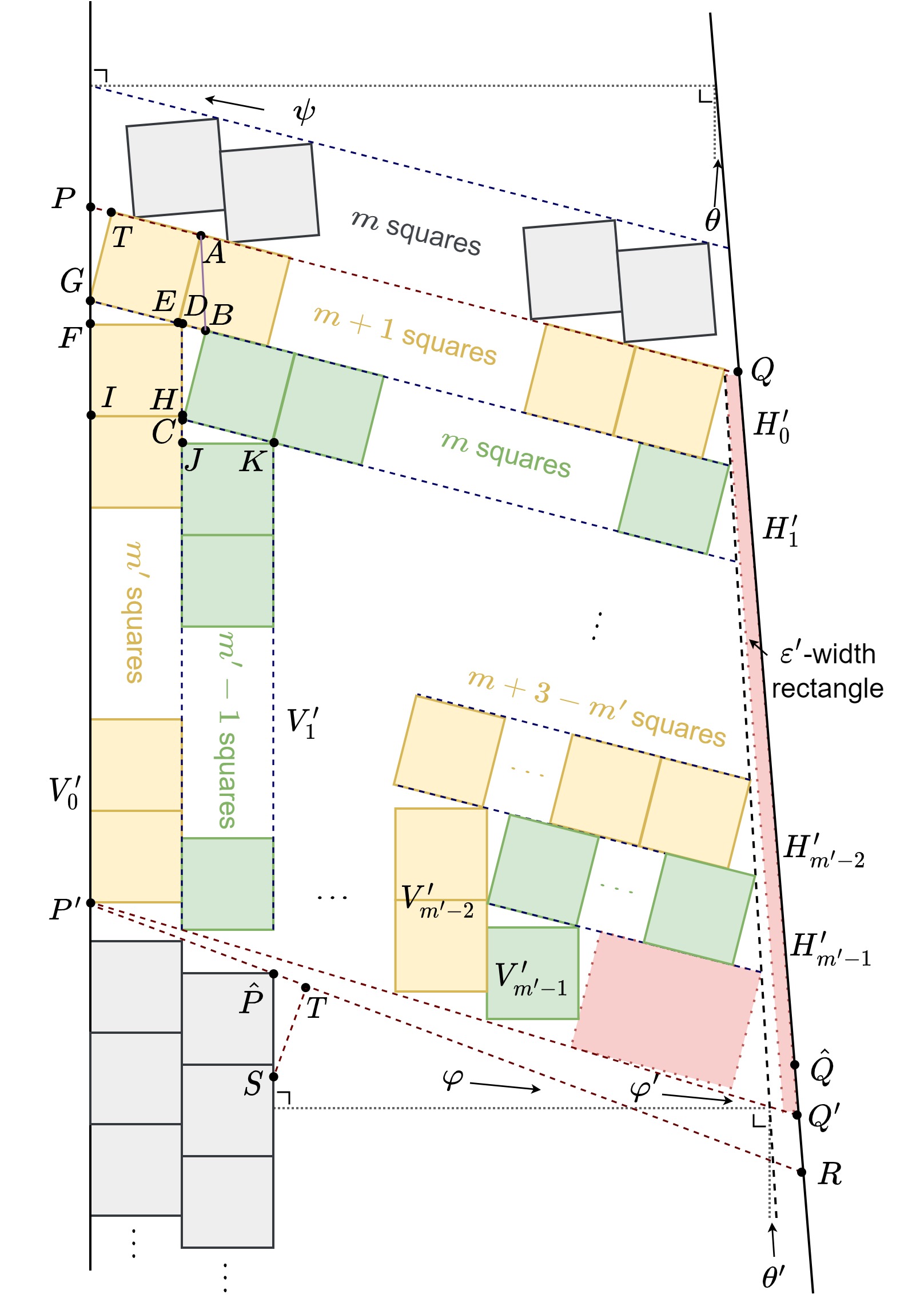}
    \centering
    \caption{The second packing algorithm.}
    \label{fig-packing_alg_2}
\end{figure}

The second packing algorithm will begin at the bottom of where the first packing algorithm left off, and proceed down the trapezoid $T$. Place a horizontal stack $H_0'$ inclined at an angle of $\psi$ up against $\overline{PQ}$ that contains $m+1$ squares and so that its left edge touches the vertical wall and its right side is some distance of $\varepsilon'$ from the inclined wall (see Figure \ref{fig-packing_alg_2}). Note that our choice of $m$ in the last section is what ensures that we have enough room. We are going to place a small rectangle of width $\varepsilon' $ along the inclined wall and will not pack that region during the second packing algorithm. Note that from (\ref{eq-tq_approx}), we have
\begin{equation}\label{eq-epsilon'_approx}
    \varepsilon' \ll \varepsilon + \theta.
\end{equation}

Next, pack a vertical stack $V_0'$ with $m'$ (an integer which is less than $m$ which we will fix later) squares against the vertical wall so that it touches $H_0'$. Choose some $\theta'$ such that we can fit another horizontal stack $H_1'$ with $m$ squares snugly against $H_0'$ so that it touches both $V_0'$ and a new inclined wall of angle $\theta'$ (we will prove the feasibility of such an arrangement shortly). Label the points on these stacks as in Figure \ref{fig-packing_alg_2}.

\subsection*{Estimating \texorpdfstring{$\theta'$}{theta'}}
\begin{figure}[t]
    \includegraphics[width=.5\textwidth]{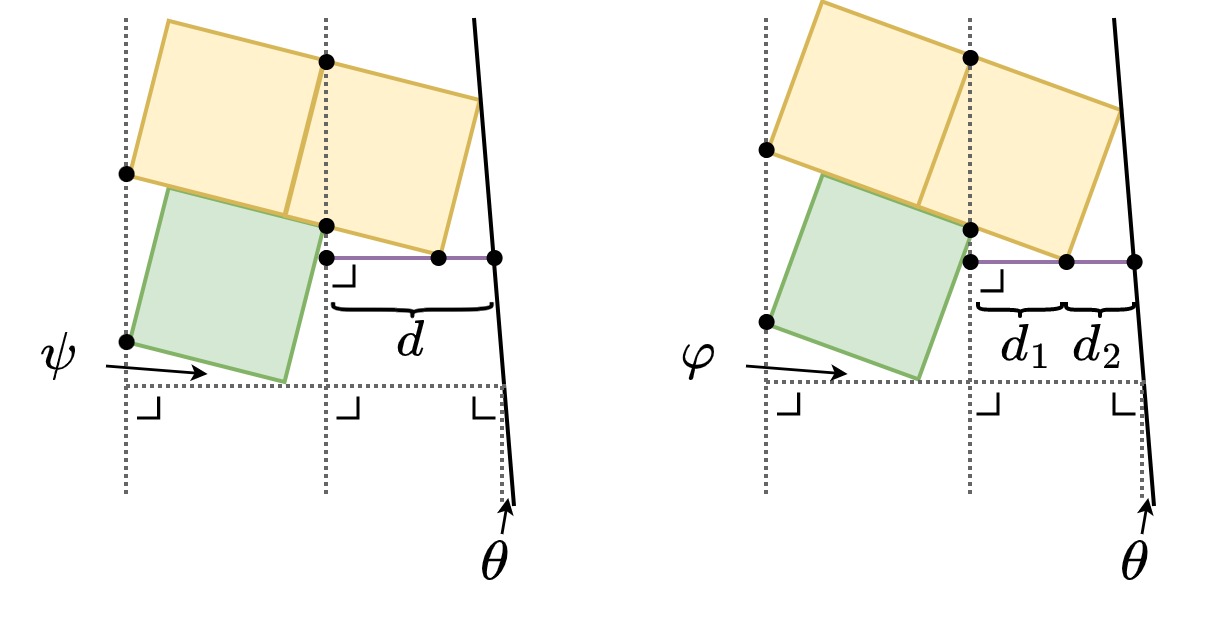}
    \centering
    \caption{Demonstrating that $\theta' < \theta$.}
    \label{fig-theta'}
\end{figure}

First we prove that $\theta' \leq \theta$ by ensuring that there is enough room to pack $H_1'$ when $\theta' = \theta$ (if the horizontal stacks were too large we would be forced to choose $\theta' < \theta$). However, this is feasible if and only if the distance $d$ depicted in Figure \ref{fig-theta'} is greater than or equal to $1$ (otherwise we would not be able to slide in the vertical stack $V_0'$). But, since $\varphi > \psi$, then $d_1 + d_2  > d$. Clearly, though, $d_1$ is greater than $DC$ in the first packing algorithm (see Figure \ref{fig-packing_alg_1}) and $d_2$ is greater than $CB$. Thus, $d > d_1 + d_2 > DB = 1$, implying that $\theta' \leq \theta$.

We now give an upper bound on the discrepancy between $\theta'$ and $\theta$.

Observe that
    \[\angle EAB = \psi + \theta', \quad \angle BCD = \psi, \quad \text{and} \quad \angle GDF = \psi. \]
Thus,
    \[ED + DB = \tan (\psi + \theta'), \quad 1 + ED= \frac{1}{\cos \psi}, \quad \text{and} \quad DB = \tan \psi.\]
We can then define $\theta'$ to be the unique solution to
\begin{equation}\label{eq-theta'}
    \sec \psi - 1 + \tan\psi = \tan (\psi + \theta').
\end{equation}
Taylor series expansion gives
    \[\theta'  =\frac{1}{2}\psi^2 + O(\psi^3) .\]
Comparing with (\ref{eq-varphi_approx}) and (\ref{eq-psi_approx}) gives
\begin{equation}\label{eq-theta'_approx}
    0 \leq \theta-\theta' \ll \varphi^3.
\end{equation}

We then can continue to pack the $m+1 - i$-square stacks $H_i$ for $i=0, 1, \dots, m'-1$ in the same manner. Similarly, we will pack the $m'-i$-square stacks $V_i$ for $i=0, 1, \dots, m'-1$. 

\subsection*{Estimating \texorpdfstring{$\varphi'$}{phi'}}
To determine the angle $\varphi'$, note that $\varphi' = \angle JIH$ (see Figure \ref{fig-packing_alg_2}). However the triangle $CKJ$ is congruent to the triangle $CDB$, which implies that $CJ = DB$. Furthermore, the triangle $CDB$ is congruent to the triangle $FGD$, which implies that $ED = HC$. Thus the triangle $JIH$ is congruent to the triangle $BAE$, which then implies that $\varphi' = \psi + \theta'$. From (\ref{eq-psi_approx}) and (\ref{eq-theta'_approx}), this implies that
\begin{equation}\label{eq-varphi'_approx}
    \varphi - \varphi' \ll \varphi^3.
\end{equation}
Note that $\varphi'$ is fixed independent of $\varepsilon$.

\subsection*{Estimating \texorpdfstring{$m'$}{m'}}
We want to choose $m'$ such that if we packed a new horizontal stack of length $m$ at an angle of $\varphi$, there would still be just enough room for two vertical stacks up against the vertical wall. Note that these stacks would continue to the bottom of the trapezoid $T$. This setup would allow us to once again apply the first packing algorithm.

Label the points as in the bottom of Figure \ref{fig-packing_alg_2}. We are assuming here that $\varphi \geq \varphi'$. The case for $\varphi < \varphi'$ is analogous. Note that we have chosen points $S$ on the right-hand side of the two vertical stacks and $T$ on $\overline{PR}$ such that $ST = 1$ and $\angle \hat{P}TS$ is a right-angle (equivalently, $\angle TS\hat{P} = \varphi$). Thus, we will want to choose $m'$ to be the least integer such that $TR \geq m$. Of course, we will have to show that such an $m'$ exists and is less than or equal to $m$, and we will have to bound the discrepancy between $m$ and $m'$. For now, it is easy to see that $m' \asymp m \asymp \theta^{-1}$.

Note that $m'$ will have to depend upon $\varepsilon$, since we are attempting to close the gap caused by the $\varepsilon'$-width rectangle. This is important as the discrepancy $TR - m$ will generate the $\varepsilon$ for the next iteration of the first and second packing algorithms, and we do not want these errors aggregating over each iteration of this pair of packing algorithms.

Let $\hat{Q}$ be the point on the inclined wall such that $\overline{P'\hat{Q}}$ is parallel to $\overline{PQ}$ (namely, inclined at an angle of $\psi$). From the sine rule, and using the fact that $m', m \ll \theta^{-1}$, we have
    \[P'\hat{Q} - PQ = PP' \frac{\sin \theta}{\sin (\pi / 2 - \psi - \theta)} = (m' + O(1)) \frac{\sin \theta}{\cos (\psi + \theta)} = m' \theta + O(\theta).\]
Again, by the sine rule,
    \[P' R   = P'\hat{Q} \left( \frac{\sin (\pi/2 + \theta + \psi)}{\sin (\pi /2 - \theta - \varphi)} \right) = \left( PQ + m' \theta + O(\theta) \right) \left( \frac{\cos (\theta + \psi)}{\cos ( \theta + \varphi)} \right)\]
Observe that since $\varphi - \psi = \theta + O(\varphi^3)$ from (\ref{eq-psi_approx}), we have
    \[\frac{\cos (\theta + \psi)}{\cos (\theta + \varphi)} = \frac{\cos(\theta + \varphi) \cos (\theta + O(\varphi^3)) + \sin (\theta + \varphi) \sin(\theta + O(\varphi^3)) }{\cos (\theta + \varphi)} = 1 + \theta \varphi + O(\varphi^4).\]
Thus, since $m' \ll m \ll \varphi^{-2}$ and $PQ = m + \varphi + O(\varphi^2 + \varepsilon)$ (from (\ref{eq-tq_approx})), we have
    \[P' R = \left(m +1 + \varphi +m' \theta + O(\varphi^2 + \varepsilon)  \right) \left( 1 + \theta \varphi + O(\varphi^4) \right) = m + 1 + \varphi + m'\theta  + \theta \varphi m + O(\varphi^2 + \varepsilon) \]
We then have
\begin{equation}\label{eq-tr_approx}
    TR = P'R - 2 (1 - \cos \varphi) - \sin (\varphi) = m + m' \theta + \theta \varphi m - 1 + O(\varphi^2 + \varepsilon). 
\end{equation}
Now, observe that it is always possible to find a $m' \leq m$ such that $TR \geq m$ as long as $\varepsilon \ll \varphi^2$, since taking $m' = m $ would make $TR \sim m+ \varphi$ (larger than $m$), so such an $m'$ must exist. We are implicitly using the fact that we can always fix $TR$ to the nearest $\theta$ (the slope of the inclined wall). We can summarize these observations with the following equation
\begin{equation}\label{eq-tr}
    0 \leq TR -m \ll \theta,
\end{equation}
noting that the asymptotic is independent of $\varepsilon$. Substituting this into (\ref{eq-tr_approx}), and once again assuming that $\varepsilon = O(\theta)$, gives
    \[\theta m' + \theta \varphi m - 1 = O(\theta).\]
This allows us to solve for $m'$:
    \[m' = \theta^{-1} + O(\varphi m ).\]
Combining this with (\ref{eq-m_approx}) gives
\begin{equation}\label{eq-m'}
    m - m' \ll \varphi^{-1}.
\end{equation}

\subsection*{Summary} We have now developed two packing algorithms. When applied back to back, the end of the second packing algorithm is the setup for the initiation of the first packing algorithm. Note that the only parameter that will vary throughout this iterative process is $m'$, as it is dependent upon $\varepsilon$.

\section{Proof of Theorem \ref{thm-main}}\label{s-proof}
\subsection*{Packing \texorpdfstring{$T$}{T}}
\begin{figure}[t]
    \includegraphics[width=.8\textwidth]{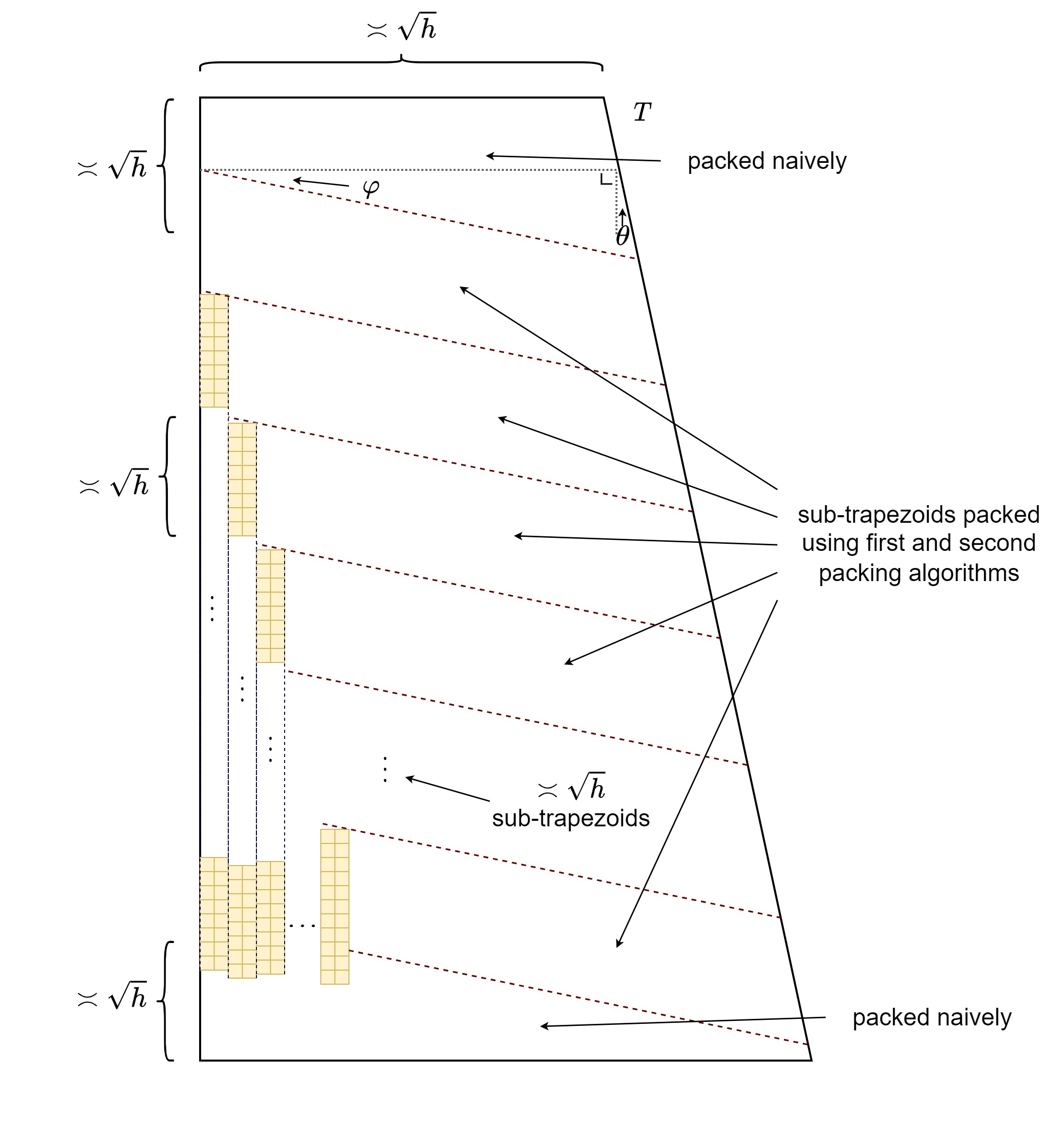}
    \centering
    \caption{Packing the trapezoid $T$.}
    \label{fig-t_packing}
\end{figure}

We can now apply our two packing algorithms to complete the proof of Theorem \ref{thm-main}. Recall that we had not yet fixed the width of the trapezoid $T$. Now, the inclined wall of $T$ was inclined at an angle $\theta \asymp h^{-1/2}$ (see (\ref{eq-theta_approx})). Observe that with this choice of $\theta$, we must have $m,m'\asymp h^{1/2}$, and $\varphi,\psi,\varphi' \asymp h^{-1/4}$. We now fix the width of the upper edge of $T$ so that the first $m$-length horizontal stack inclined at an angle of $\varphi$ is within $O(\sqrt{h})$ of the top of $T$, thereby implying that $T$ has a width on the order of $\sqrt{h}$ (see Figure \ref{fig-t_packing}).

To pack $T$, we apply our two packing algorithms back to back, packing a "sub-trapezoid" at the top of $T$ of height about $2\theta^{-1} \asymp \sqrt{h}$. Note that our first choice of $\varepsilon$ can be $0$. We then pack the two vertical stacks along the \textit{entire} vertical wall. We repeat this process, packing each sub-trapezoid all the way until we have packed everything except a region at the bottom of $T$ that has height at most $O(\sqrt{h})$ (see Figure \ref{fig-t_packing}). Note that at each stage, we are using the same $\varphi, \psi, \varphi'$ and $m$, but that the choices for $\varepsilon, \varepsilon'$ and $m'$ will change for each iteration of the packing algorithm. However, we can always ensure that $\varepsilon \ll \theta$, courtesy of (\ref{eq-tr}). Note that we will have to apply our algorithm a total of $\asymp \sqrt{h}$ times to pack all of $T$.

\subsection*{Estimating the wasted space}
Let us begin by computing the wasted space from a single iteration of the first and second packing algorithms. Observe that \textit{all} of the rectangles in the wasted space have area at most $O(\frac{1}{\sqrt{\theta}})$. To see this, first observe that $H'_{m-1}$ in Figure \ref{fig-packing_alg_2} has length at most $O(\varphi^{-1})$ from (\ref{eq-m'}). Second, note that we can always force $\varepsilon \ll \theta$, and thus $\varepsilon' \ll \theta$ (see (\ref{eq-tr}) and (\ref{eq-epsilon'_approx})).

Next, note that there are two large sliver triangles of length $O(\theta^{-1})$ and angles $\theta - \theta'$ and $\varphi - \varphi'$ (see Figure \ref{fig-packing_alg_2}). Thus, from (\ref{eq-theta'_approx}) and (\ref{eq-varphi'_approx}), the area of these triangles is $\ll \frac{1}{\sqrt{\theta}}$. There is also a small sliver between $\overline{P'Q'}$ and the rectangle of width $\ll \varphi^{-1}$. This has an angle trivially bounded by $\varphi$, and so its area is $O(\frac{1}{\sqrt{\theta}})$.

Finally, observe that there are $\ll \theta^{-1}$ $O(1)$-length triangles created through both packing algorithms, and their angles are all $\ll \sqrt{\theta}$, implying that there total contribution is $\ll \frac{1}{\sqrt{\theta}}$.

Thus, the total amount of wasted space generated through these two packing stages is $\ll \frac{1}{\sqrt{\theta}} \ll h^{1/4}$. Since we are applying these packing algorithms $\asymp \sqrt{h}$ times, then the total wasted space contributed by packing the sub-trapezoid regions becomes $O(h^{3/4})$. We pack the regions at the top and bottom of $T$ trivially, which generates an additional wasted space of $O(h^{1/2})$ (see Figure \ref{fig-t_packing}). Thus, the total wasted space generated by packing $T$ is $O(h^{3/4})$.

Now, the wasted space generated by packing the $O(x)$ vertical stacks inclined at an angle $\theta$ in Figure \ref{fig-r_packing} is $O(x \theta)$, meaning that our total wasted space is bounded as follows:
    \[W(x) \ll \frac{x}{\sqrt{h}} + h^{3/4}.\]
Equalizing the two terms gives a choice of $h \sim x^{4/5}$, yielding the desired result.

\bibliography{references}

\end{document}